\documentclass{amsart}
\usepackage{amssymb}
\usepackage{url}
\usepackage{amscd}
\newtheorem{thm}{Theorem}[section]

\newtheorem{cor}[thm]{Corollary}
\newtheorem{lemma}[thm]{Lemma}
\newtheorem{prop}[thm]{Proposition}
\theoremstyle{remark}

\theoremstyle{definition}

\newtheorem{ex}[thm]{Example}
\numberwithin{equation}{section}
\numberwithin{thm}{section}
\def\bn#1#2{B_{#2}({o_{#1}})}
\def\R{\text{$\mathbb{R}$}}
\def\fis#1{\dot{#1}}
\def\lra{\longrightarrow}
\def\pt#1#2{\tau^{#2}_{#1}}
\def\nb{\nabla}

\def\nbp#1#2{\nabla_{\frac{\partial}{\partial #1}}#2}
\def\tpm{T_p M}

\def\pt#1#2{\tau^{#2}_{#1}}
\def\e#1{\exp_{#1}}

\def\bn#1#2{B_{#2}({0_{#1}})}
\def\p1#1#2{\frac{\partial #1}{\partial #2}}
\def\to{t_o}
\def\pr{\parallel}
\def\cal#1{\mathcal{#1}}
\def\ov{\overline}
\def\taT{\tilde{T}}

\def\ptn#1#2{\chi^{#2}_{#1}}

\def\D{D^2 E(c)}
\def\so{\rightarrow}
\newcommand{\ed}{\end{document}}

\newcommand{\We}{A_{\xi}}
\def\Ex{\mathrm{Exp}^{\perp}}

\def\N{\text{$\mathbb{N}$}}
\def\E{\text{$\mathbb{E}$}}
\def\H{\text{$\mathbb{H}$}}
\def\pr{\parallel}
\def\d1#1#2{\frac{d#1}{d#2}}
\begin{document}
\title[Weak Riemannian geometry]{Exponential map of
a weak Riemannian Hilbert manifold}
\author{Leonardo Biliotti}
\address{Via Turati 10 \\ 50136 Firenze Italia}
\email{biliotti@math.unifi.it }
\thanks{Research partially supported by CNPq  (Brazil)}
%
%
\begin{abstract}
We prove the Focal Index Lemma and the Rauch and Berger comparison Theorems
on a weak Riemannian  Hilbert manifold with a smooth Levi-Civita
connection and we apply these results
to the free loop space $\Omega (M^n)$ with  the  $L^2$
(weak) Riemannian structure.
\end{abstract}
\maketitle
\section{Introduction}

As a preliminary step to understand the global geometry of
a Riemannian Hilbert manifold $M$, one studies singularities of
its exponential map. Singular values of $\exp$ are the conjugate
points in $M.$ In infinite dimension, 
there exist two types of conjugate points: when
the differential of the exponential map fails to be injective
(a monoconjugate point) or when the differential of the exponential
map fails to be surjective (a epiconjugate point). More ge\-ne\-ral\-ly, let 
$N$ be a submanifold of $M$ such that $\forall p \in N$ the 
tangent space at $p$ of $N,$ $T_p N,$ is a closed subspace of $\tpm.$ 
Singular values of the map 
$\Ex : T^{\perp} N \lra M,$ defined by 
$\Ex (X)= \exp (X),$ where $\exp$ 
is the exponential map of $M$ and
$T^{\perp} N$ is the normal bundle of $N$, are called focal points: 
a \emph{monofocal point}, 
when the differential fails to be injective, and an \emph{epifocal point}, 
when the
differential fails to be surjective. Clearly, the application 
$\Ex$ is defined a priori 
only  in an open subset which contains the ``zero section'', 
i.e. 
the subset $ \{ 0_p \in T_p^{\perp} N:\ p \in N \} \subseteq T^{\perp} N$. 

Let now $(M,\langle\cdot ,\cdot \rangle) $ be a \emph{weak} Riemannian
Hilbert manifold with a smooth Levi-Civita connection $\nabla$,
whose existence is not guaranteed a priori. It defines  parallel
transport, curvature tensor $R$, geodesics and a smooth
exponential map. 
These manifolds have
been intensively studied and they have found many diverse applications
particularly in geometry, calculus of variations 
and mathematical physics 
(see \cite{arnold},
\cite{em}, \cite{df}, \cite{mi2}, \cite{Mi3}, 
\cite{Mi2}, \cite{Mi4}, \cite{ps} ). For
example, see \cite{arnold}, \cite{em}, \cite{Mi3}, \cite{Mi4}, 
any motion of a perfect
fluid cor\-responds to a geodesic on the group of volume-preserving 
diffeomorphism of a compact manifold $M,$ which is the region filled with 
fluid, with respect to the weak Riemannian me\-tric
which is given by the $L^2$ inner product on each tangent space. Moreover,
existence of conjugate points are   
related to stability of the fluid flows of $M$.

An other  important example is 
the free loop space $\Omega (M^n )$ of a compact
manifold $M^n,$ which is among simplest Hilbert manifold.
This Hilbert manifold has been intensively studied and it has many diverse
applications (see \cite{df}, \cite{Mi2}, \cite{ps}). It has an
$L^2$ metric that is a weak Riemannian structure, which induces
a smooth Levi-Civita connection and a smooth exponential map. 

One motivation for the results
presented here was the paper of Misiolek \cite{Mi2}  
where it was proved that for every $s>0$, 
the exponential map of the $H^s$ metric on the 
$H^{s_o} (S^1, G),$ i.e. the set of Sobolev $H^{s_o}$
maps from the unit circle $S^1$ into a compact, connected Lie group $G,$
is a nonlinear Fredholm map of zero index while the exponential map of 
$\Omega ({\rm SU}(2))$
with respect to the $L^2$ weak Riemannanian metric is not.  

Since the model space $\H,$ on which $M$ is modeled,
is a Hilbert space,  it is
possible to transport on the tangent space the structure of a
topological vector space, that we will denote by $\tau,$ given by
the chart and this topology can be induced by a scalar product
(see \cite{La} page 26). We assume that the curvature tensor $R$
is a trilinear continuous operator with respect to the topology
$\tau.$

Let $N$ be a submanifold of $M$ such that for some $p \in N,$ $T_p
N$ is a closed subspace of $(\tpm, \tau )$ and $T_p M = T_p N
\oplus T_p^{\perp} N.$ In this context, we shall define the notion
of focal point along a normal geodesic starting from $p,$ which is
equivalent to the usual one in Riemannian geometry. We shall
prove the Focal Index Lemma, when there exist a finite number of epifocal 
points which are not
monofocal along a geodesic of finite length, which generalizes  the Index
Lemmas, see \cite{CE} page 24, in finite dimensional Riemannian
geometry. As immediate corollaries we 
get the Rauch and Berger comparison Theorems.

After formulating and proving the Focal Index Lemma, and its
corollaries, we apply it, in Section $4,$ to the loop group $\Omega (M^n).$ 
We prove that a geodesic $c:[0,b] \lra \Omega (M^n)$ 
with length big enough has
conjugate points and its index is infinite. A similar result can be proved 
for focal points of $c(0)$ along $c$ with respect to the geodesic submanifold
defined by $\fis{c}(0)$. Then we
analyze the case when $M^n =G$  is a non-abelian compact Lie group and we
prove that 
its exponential map fails to be Fredholm. Moreover,  we give an
example of a submanifold $N$ of $\Omega (G)$ such that 
$\Ex : T^{\perp} N \lra M,$ fails to be Fredholm as well.
\section{Exponential map on  Hilbert Manifolds}
\mbox{} In this section we will recall some general results and
well known facts. Our basic references are \cite{AMR}, \cite{Kl}
and \cite{La}. 

Let  $M$ be a Hilbert manifold modeled on an
infinite dimensional Hilbert space $\H.$ Recall that a weak
Riemannian metric on $M$ is a smooth assignment to each point $p
\in M$ of a continuous, positive definite, symmetric bilinear form
$p \lra \langle \cdot, \cdot \rangle (p)$ on the tangent space
$\tpm.$ Note that $\tpm \cong \H$ need not be complete as a metric
space under the distance induced by $\langle \cdot, \cdot\rangle
(p).$ Consequently the existence of a smooth Levi-Civita
connection $\nabla$ associated with a weak Riemannian metric
is not immediately guaranteed. If, however, such a connection
exist, it is necessarily unique.

Throughout this paper we shall assume that $M$ is a Hilbert
manifold endowed with a weak Riemannian metric $\langle \cdot
,\cdot \rangle, $ and $M$ will be called weak Riemannian Hilbert
manifold. We further assume that $M$ admits a Levi-Civita
connection $\nabla$, whose curvature tensor $R$ is a continuous
trilinear operator of the tangent space with respect to the
topology $\tau.$

For any $p \in M$ the exponential map $\e p : \tpm \lra M$ is a local
diffeomorphism in a neighborhood of the origin in $\tpm.$ The differential
$ d (\e p ) $ can be computed using the Jacobi equation that is
the linearized version of the geodesic equation.

Let $c :[0,b] \lra M$ be a geodesic.
A vector field $J$ along $c$ is called Jacobi
field if it satisfies the Jacobi differential equation
$$
\nbp{t} \nbp{t}{J}(t) + R(J(t), \fis c (t) ) \fis c (t) =0,
$$
where $\nbp{t}{}$ denotes the covariant derivation along $c.$

It is well known that if $c (t) = \e p (tv)$ is the geodesic starting at $p$
in the direction $v,$ then the vector field
$Y(t) = d (\e p )_{tv} (tw) $ satisfies the Jacobi differential
equations  with  initial values $Y(0)=0$ and
$\nbp{t}{Y} (0)=w.$

Let $N$ be a subma\-ni\-fold of $M$ and let $c:[0,b] \lra M$ be a geodesic
such that $c(0)=p \in N$
and $\xi=\fis{c} (0) \in T_p^{\perp} N,$ i.e. $c$ is a normal 
geodesic of $N$. Suppose also
that $T_p N$ is a closed subspace of $(\tpm, \tau)$ and $\tpm =
T_p N \oplus T_p^{\perp}N.$ This happens when $N$ is a submanifold
of $M$ defined by some vector $v \in \tpm$: $N=\e p
(\bn{p}{\epsilon} \cap <v>^{\perp}),$ where $\epsilon$ is
sufficiently small such that the map $\e p : \bn{p}{\epsilon} \lra M$ is a
diffeomorphism onto the image. 

As in the Riemannian case, the Weingarten operator  
is given by $\We (X)= -P (\nabla_{X} \xi (p) ),$ 
where $P$ is the projection of $\tpm$ onto $T_p N.$
Of course, the Weingarten operator is  a linear
continuous map of $(\tpm , \tau)$ and symmetric with respect to
$\langle \cdot, \cdot \rangle.$ In finite dimensional Riemannian geometry, the 
Jacobi fields  
along $c$ with initial values 
\[
J(0) \in T_p N, \   
\nbp{t}{} J (0) + \We (J(0)) \in T_p^{\perp} N,
\] 
which are called 
$N-$Jacobi fields, describe completely the differential of the map 
$\Ex .$   

Now, let $ \pt{t}{s} : T_{c(t)} M
\lra T_{c(s)} M$ be the isomorphism between the tangent
spaces given by parallel transport along a geodesic $c.$
Since the parallel transport along $c$ and $\nbp{t}{}$ commute, we can
rewrite the Jacobi equation relative to $N,$
as an initial value problem on $\tpm$ as follows
\[
\left \{ \begin{array}{l}
T''(t)\ +\ R_t(T(t))=0;\\
T(0)(v,w)=(v,0),\ T'(0)(v,w)=(-\We (v),w),
\end{array}
\right.
\]
where
$$
R_t:\tpm    \lra \tpm, \ \
R_t(X)=\pt{t}{0}(R(\pt{0}{t}(X),\fis{c}(t))\fis{c}(t))
$$
is a family of symmetric operators of $\tpm$.
We will call the above differential equation {\em Jacobi flow } of
$c$ relative to $N.$
One may also note that all maps $\Phi(t)$ defined by
$$
\begin{array}{lcl}
& \tpm \times \tpm \stackrel{\Phi(t)}{\lra} \R & \\
& (u,v) \lra  \langle T(t)(u),T'(t)(v) \rangle & \\
\end{array}
$$
are symmetric; indeed $\Phi (0)$ is symmetric since $\We$ is a
symmetric operator, and
\[
( \langle T(t)(u),T'(t)(w) \rangle \ - \  \langle T(t)(w),T'(t)(u)
\rangle )' =0.
\]

A point $q=c(t_o)$ is called a {\em monofocal point} 
respectively an {\em epifocal point} of $p=c(0)$ along $c$ 
if $T(t_o)$ fails to be injective 
respectively fails to be surjective. In general, we call a point $q=c(t_o)$ 
a \emph{focal point} of $p$ along $c$ when $T(t_o)$ is not an isomorphism. 
One may note that this definition is equivalent to the
Riemannian one. 

Now, let $\E_1$ and $\E_2$ be Banach spaces. A bounded linear
operator  $T:\E_1 \lra \E_2$ is called Fredholm if it has a closed
range and its kernel and co-kernel ($coker T= \E_2 / T( \E_1))$ are finite
dimensional. The index of $T$ is the number ${\rm ind} \ T = \dim
Ker T - \dim coker T.$

A smooth map between Banach manifolds $f:M \lra S$ is called
Fredholm if for each $p \in M$ the derivative $d(f)_p: T_p M \lra
T_{f(p)} N$ is a Fredholm operator. If $M$ is connected then the
${\rm ind}\ (df)_p$ is independent of $p$, and  one defines
the index of $f$ by setting ${\rm ind} (f)={\rm ind} (df)_p$ (see
\cite{et}, \cite{smale}). One may note that if
the map $\Ex$ is well defined, then $\Ex$ is a nonlinear Fredholm 
map if the Jacobi flow along any normal geodesic of $N,$ 
describes a curve in the Fredholm operators for every $t>0.$

Now, we shall describe  the adjoint operator of $T(b),$ since we shall
understand the behavior of the focal points of $c(0)$ along the geodesic
$c$.

Let  $u \in \tpm$ and let $J$ be the Jacobi field along the geodesic
$c$ such that
$
J(b)=0, \ \nbp{t}{J(b)}=\pt{0}{b}(u).
$
By a lemma from Ambrose, see \cite{Am} or \cite{La} Lemma 3.4 page 243, 
we have
\[
(1)\ \langle T(b)(v,w),u \rangle = \langle T(0)(v,w),
\nbp{t}{J}(0) \rangle \ - \ \langle
T'(0)(v,w),J(0) \rangle.
\]
Let $\ov{c}(t)=c(b-t).$ Let
\[
\left \{ \begin{array}{l}
\taT''(t)\ +\ R_t (\taT(t))=0;\\
\taT(0)=0,\ \taT'(0)=id,
\end{array}
\right.
\]
be the Jacobi flow of $\ov{c}$ relative to the submanifold 
$\ov{N}=\{ \ov{c} (0) \}$.
It is easy to check that if $J$ is a Jacobi field along $c,$ then
$\overline{J}(t)=J(b-t)$ is the Jacobi field along
$\ov{c}$ such that
$\nbp{t}{\overline{J}}(b)=-\nbp{t}{J}(0)$. Then $(1)$ becomes
$$
\begin{array}{lcl}
\langle   T(b)(v,w),u \rangle &=&  \langle  (v,0),
\pt{b}{0}(\taT'(b)(-\pt{0}{b}(u))) \rangle \\
&-&  \langle (- \We (v),w),\pt{b}{0}(\taT (b)(-\pt{0}{b}(u))) \rangle, \\
\end{array}
$$
so the adjoint operator is given by
$$
\begin{array}{lcl}
\langle T^*(b)(u),(v,0) \rangle  &=& -   \langle \pt{b}{0}
 ( \taT'(b) (\pt{0}{b}(u) ) ) , (v,0) \rangle \\
 & + &
\langle \We ( P ( \pt{b}{0}(\taT(b) ( \pt{0}{b}(u) ) ) ) ) , (v,0) \rangle,
\\
\langle  T^*(b)(u), (0,w) \rangle  & = &  \
\langle \pt{b}{0}(\taT(b)(\pt{0}{b}(u))), (0,w) \rangle .    \\
\end{array}
$$
\begin{prop} \label{p1}
The kernel of $ T(b)$ and the kernel of $T^* (b)$ are isomorphic.
\end{prop}
\begin{proof}
$\ $ Let $w \in \tpm$ be such that $T(b)(w)=0.$ The Jacobi field 
$Y(t)=\pt{0}{t}(T(t)(w))$ vanishes at $t=b.$ Then there exists a unique 
$w \in  T_{c(b)} M$ such that
$$
Y(b-t)=\pt{b}{b-t}(\taT(t)(\overline{w})).
$$
Using the boundary conditions
of $Y(t)$ we get $T^* (b) (\pt{b}{0} (\ov{w}))=0.$ 
In particular the application
$$
\begin{array}{lcl}
& f_1: Ker T(b) \lra Ker T^* (b) & \\
& w \lra \pt{b}{0} (\ov w ) & \\
\end{array}
$$
is an injective linear map. 

Vice-versa, let $v \in Ker T^* (b)$. 
We denote by $\ov v = \pt{0}{b} (v)$ and we consider the following 
Jacobi field 
$Y(t)=\pt{b}{b-t} ( \tilde{T} (t) (\ov v ))$ along $\ov{c}$. 
Since $T^* (b) (v)=0,$ there exists a unique $\theta \in \tpm$ such that
$Y(b-t)=\pt{0}{t} (T(t)(\theta))$. Hence $T(b) (\theta)=0$ since $Y(0)=0;$
moreover, the application
$$
\begin{array}{lcl}
& f_2: Ker T^* (b) \lra Ker T (b) & \\
& v \lra \theta & \\
\end{array}
$$   
is injective and  one may note that $f_2  \circ f_1 =Id,$ 
thus concluding our
proof. 
\end{proof}
\begin{prop} \label{singular1}
Let $(M,\langle \cdot,\cdot \rangle )$ be a weak Riemannian
Hilbert manifold with a smooth Levi-Civita connection $\nabla$. Let $N$ be
a submanifold of $M$ and let $c:[0,b] \lra M$ be a normal geodesic of $N$, i.e.
$c(0)=p \in N$ and $\xi=\fis{c} (0)\in T_p^{\perp} N.$ 
Assume also that $T_p N$ is a closed subspace of 
$(\tpm, \tau)$ such that
$T_p M = T_p N \oplus T_p^{\perp} N.$ Then we have
\begin{enumerate}
\item if $c (t_o)$ is not a monofocal of $p$ along $c$, 
then the image of $T(t_o)$ is 
a dense subspace relative to the topology $\tau$ induced by the metric
$\langle \cdot, \cdot \rangle(p);$ 
\item if $c(t_o)$ is a monofocal of $p$ along $c$ 
then $c(t_o)$ is an epifocal of $p$ along $c$; 
\item when $N=\{p \}$ a point $q=c (t_o)$ is
a monofocal of $p$ along $c$ if and only if $p$ is a monofocal point 
of $q$ along  $\ov{c}(t)=c(t_o-t);$ 
\item when $N=\{p \}$  if a point $q=c (t_o)$
is an epifocal of $p$ along $c$ 
and the image of $T(\to) $ 
is a closed subspace of $(\tpm, \tau)$, then $p$ is a monofocal point 
of $q$ along
$\ov{c}(t)=c(t_o-t);$
\end{enumerate}
\end{prop}
\begin{proof} $\\ $
\begin{enumerate}
\item since $Ker T^* (t_o)=0,$ by Proposition \ref{p1},
the closure of $Im T(t_o)$ with respect to $\tau'$ satisfying 
$\ov{Im T (t_o) }^{\perp}=0.$ Now, one may note that $\tau'$ makes $\tpm$
into a locally convex space and applying
3.10 and 3.5  Theorems in \cite{rudin} we have the statement.
\item one may note that
the adjoint of the Jacobi flow of $c$ is the
Jacobi flow of $\ov{c};$ 
\item since the image of $T (t_o)$ is closed, then $q$ is a monofocal point
of $p$ along $c$ as well. Now, the statement follows from the above item.
\end{enumerate}
\end{proof}
\section{Focal Index lemma and  Rauch and Berger comparison
theorems in weak Riemannian geometry}
\mbox{}
In the infinite dimensional case,
the distribution of singular points  of the exponential map along a
geodesic of finite length is different from the finite dimensional case.
Indeed, Grossman showed \cite{Gr} how the distribution
of monoconjugate points can have cluster points.
The following example proves that the same situation
may oc\-cur in the case of focal points along  a
geodesic of finite length.
\begin{ex} \label{esempio}
Let
$M= \{ x \in l_2 :\ x^2_1\ +\ x^2_2 \ + \ \sum_{i=3}^{\infty} a_i x^2_i=1 \}$,
where $( a_i )_{i \in \N}$ is a positive sequence of real numbers.
$M$ is a Riemannian Hilbert manifold and it is easy to check that
\[
\gamma(s)= \sin (s) e_1 \ + \ \cos (s) e_2
\]
is a geodesic and $T_{\gamma(s)}M= < \fis{\gamma}(s),e_3,e_4,\ldots
>$. Let $N$ be a subma\-ni\-fold defined by $\fis{\gamma}
(0).$ We shall restrict ourselves to the normal $N-$Jacobi fields, i.e. the 
Jacobi fields which satisfy $\langle J(0), \fis{c}(0) \rangle=0$ . Since
for $k \geq 3$
\[
E_k:= \{ x^2_1\ +\ x^2_2 \ +\ a_k x^2_k =1\ \} \hookrightarrow M
\]
is a closed totally geodesic submanifold of $M$, the 
sectional curvature of the plane
$<\fis{\gamma}(s),e_k>$ is given by $K(\fis{\gamma}(s),e_k)=a_k$ and 
consequently the Jacobi fields with boundary conditions $J_k(0)=e_k$, 
$\nbp{t}{J_k(0)}=0$, are given by $J_k(t)=\cos (\sqrt{a_k}t)e_k$.
Hence
\[
d (\Ex)_{s \fis{\gamma}(0)} (\sum_{k=3}^{\infty} b_k e_k )=
\sum_{k=3}^{\infty} b_k \cos (\sqrt{a_k}s) e_k.
\]
Clearly, the points $\gamma (r_k^m)$, where $r_k^m=\frac{m \pi}{2 \sqrt{a_k}},$
$m \in \N,$
are monofocal of $e_2$ along $\gamma$. Specifically, let
 $a_k=(1- \frac{1}{k})^2.$ The points $\gamma (s_k)$,
where $s_k=\frac{k \pi}{2 (k-1)},$ are monofocal of $e_2$ along 
$\gamma$, $s_k \so
\frac{\pi}{2}$ and
\[
d (\Ex)_{ \frac{\pi}{2} \fis{\gamma}(0)} ( \sum_{k=3}^{\infty} b_k
e_k)=
      \sum_{k=3}^{\infty} b_k  \cos (\frac{k-1}{k} \frac{\pi}{2}) e_k.
\]
Hence $\gamma(\frac{\pi}{2})$ is not monofocal of $e_2$ along $\gamma.$ On
the other hand if $\sum_{k=3}^{\infty} \frac{1}{k} e_k = d
(\Ex)_{\frac{\pi}{2}\fis{\gamma}(0)}(\sum_{k=3}^{\infty} b_k e_k)$, then $
\sin (\frac{\pi}{2k}) b_k= \frac{1}{k},$ so we have
\[
\lim_{k \so \infty} b_k = \lim_{k \so \infty} \frac{\pi}{2k}
\frac{1}{\sin(\frac{\pi}{2k})} \frac{2}{\pi}=\frac{2}{\pi}.
\]
This means that  $\gamma(\frac{\pi}{2})$ is  an epifocal point 
of $e_2$ along $\gamma$.
\end{ex}
This example shows that there exist  epifocal  points which are not
monofocal. We call them
\emph{pathological points}.
Clearly, if  the exponential map is a non-linear Fredholm map,
and then it is necessarily of zero index,
monoconjugate points and epiconjugate points
along geodesics  coincide. This holds for
the Hilbert manifold $\Omega (M^n)$ 
(see \cite{Mi}) with the $H^1$ Riemannian structure.

Now we shall prove  the Focal Index Lemma.

Let $N$ be  a
submanifold of $M$ and let $c:[0,b] \lra M$ be a geodesic of $M.$
Assume that
$c(0)=p \in N,$ $\xi=\fis{c} (0)\in T_p^{\perp} N,$
$T_p N$ is a closed subspace of $(\tpm , \tau)$ and finally
$T_p M = T_p N \oplus T_p^{\perp} N.$

Let $X:[0,b] \lra \tpm$, such that $X(0) \in T_p N.$ We define
the {\em focal index form} of  $X$ as follows:
$$
\begin{array}{lcl}
I^N (X,X)& =& \int_{0}^{b}
\langle  \fis X (t), \fis X (t) \rangle - \langle R_t(X(t)),X(t)
\rangle dt
\\
& - &  \langle \We (X(0)),X(0) \rangle .
\end{array}
$$
One can note that any vector field along
$c$ is the parallel transport of  unique application
$X:[0,b] \lra \tpm.$ We will denote by
$\overline X (t)=\pt{0}{t}(X)$ the vector field along
$c$ starting from $X$.
\begin{lemma}
$I^N (X,X)= D^2 E(c)(\overline X, \overline X)$,
where  $D^2 E(c)$ is the
index form  of  $B=N \times M \hookrightarrow M \times M.$
\end{lemma}
\begin{proof} We recall that
\begin{small}
\begin{eqnarray}
\D( \overline X , \overline X ) &=& \int_0^b
\pr \nbp{t}{\overline X (t)} \pr^2
                    \ - \ \langle \overline X (t),
R(\overline X (t),\fis c(t)) \fis c(t) \rangle dt
\nonumber \\
             &-& \langle \langle A_{(\fis c(0),-\fis c(b))}(\overline X(0),
\overline X (b)),(\overline X (0),\overline X (b)) \rangle \rangle,
             \nonumber
\end{eqnarray}
\end{small}
see \cite{Sa}, where $A$ is the Weingarten operator of
$N \times M \hookrightarrow M \times M$ and 
$\langle \langle \cdot ,\cdot \rangle \rangle$ 
is the natural weak
Riemannian structure on $M \times M$ 
induced by $\langle \cdot, \cdot \rangle .$
Hence, it is enough to prove that
$\nbp{t}{\ov{X}}(t)=\pt{0}{t}(\fis X (t))$. Let $Z(t)$ be a parallel transport
of a vector $Z \in \tpm.$ Then
$$
\begin{array}{ccl}
\langle \nbp{t}{\overline X (t)},Z(t) \rangle  &= &
\langle \overline X (t),Z(t) \rangle ' \\
                          &= & \langle \fis{X} (t), Z \rangle \\
                          &= & \langle \pt{0}{t}(\fis{X}(t)), Z(t) \rangle.
\end{array}$$
\end{proof}
\begin{lemma} \label{fp} {\bf (Focal Index Lemma)}
Let $c:[0,b] \lra M$ be a geodesic with a finite number of pathological
points on its interior. Then for every vector field $Z$ along  $c$ with
$Z(0) \in T_p N,$ the index form of $X$ relative to  the submanifold
$N \times M \hookrightarrow M \times M$,
satisfies $\D (Z,Z) \geq \D(J,J)$, where
$J$  is the $N$-Jacobi field such that  $J(b)=X(b)$.
\end{lemma}
\begin{proof}
First of all we shall assume that there are no  focal points of $c(0)$ along
$c.$

We know that $Z= \ov{X}$ where  $X:[0,b] \lra \tpm.$
Since  $T(t)$ is invertible, 
there exists a piecewise differentiable application
$Y:[0,b] \lra \tpm$ such that
$Y(0)=X(0) \in T_p N$ and
$X(t)=T(t)(Y(t))$. Hence
\[
\fis{X}(t)=T'(t)(Y(t)) \ + \ T(t)(\fis{Y}(t))=A(t) \ + \ B(t).
\]
The focal index form of $X$ is given by
\begin{small}
$$
\begin{array}{ccl}
I^N (X,X) &=& \int_{0}^{b} \langle A(t),A(t) \rangle  \ + \ 2 \langle A(t),B(t) \rangle \ + \
               \langle B(t),B(t) \rangle  dt  \\
       &-&  \int_{0}^{b} \langle R_t (T(t)(Y(t))),T(t)(Y(t)) \rangle dt \ - \
\langle \We (X(0)),X(0) \rangle.
\end{array}$$
\end{small}
One can prove that
\begin{eqnarray}
\langle A(t),A(t) \rangle  &=& \langle T(t)(Y(t)),T'(t)(Y(t) \rangle'
\nonumber  \\
            &-& 2 \langle B(t),A(t) \rangle                 \nonumber   \\
            &+& \langle T(t)(Y(t)),R_t (T(t)(Y(t))) \rangle  \nonumber     \\
            \nonumber
\end{eqnarray}
since the bilinear form $\Phi (t)$ is symmetric. Hence,
the focal index form of $X$ is given by
\[
I^N (X,X) \  = \  \langle T(1)(u),T'(1)(u) \rangle \ + \
  \int_{0}^{b} \pr T(t)(\fis{Y}(t)) \pr^{2}dt.
\]
This proves the Focal Index Lemma in this case.
Moreover, if there are no  focal
points along $c,$ the focal index of a vector field $Z,$ with $Z(0) \in T_p N$
along $c$ is equal to the focal index of the $N-$Jacobi field $J$
along $c$ such that $Z(b)=J(b)$ if and only if $Z=J.$

Now, assume that there exists a pathological point on the interior of
$c;$ this means that the Jacobi flow is an isomorphism
for every $ t \neq t_o$ in  $(0,b)$ and when
$t= t_o,$ by Proposition \ref{singular1} (1), 
$T(t_o)$ is a linear operator whose image is a dense subspace.
Let $X:[0,b] \lra \tpm$ be a piecewise differentiable application with
$X(0) \in T_p N$.
Given $\epsilon>0$, there exist  $X_n^{\epsilon}$, $n=1,2$
such that
\[
\begin{array}{ccc}
& \pr T(\to)( X_1^{\epsilon})) -  X(\to)) \pr   \leq \frac{\epsilon}{4} & \\
& \pr T(\to)( X_2^{\epsilon})) -  \fis X(\to)) \pr \leq \frac{\epsilon}{4}.
\end{array}
\]
Choose  $Y^{\epsilon}$ such that
\[
\pr T(\to)(Y^{\epsilon}) -  T'(\to)(X_1^{\epsilon}) \pr \leq
\frac{\epsilon}{4} .
\]
Hence there exists  $\eta(\epsilon) \leq \frac{\epsilon}{2}$
such that for $ t \in (\eta(\epsilon) -  \to, \eta(\epsilon) +  \to)$
we have
\[
\begin{array}{ccc}
(1) & \pr  T(t)( X_1^{\epsilon}  +
   (t  -  \to)(X_2^{\epsilon}   -  Y^{\epsilon}))
    -   X(t)  \pr  \leq \epsilon, & \\
(2) & \pr \dfrac{d}{dt} (T(t)( X_1^{\epsilon}  +
(t  -  \to)(X_2^{\epsilon}  - Y^{\epsilon})))
  -  \fis X (t) \pr \leq \epsilon. &
\end{array}
\]
We denote by  $X^{\epsilon}$ the application
\begin{small}
\[
X^{\epsilon}(t)=
\left \{ \begin{array}{ll}
X(t) & {\rm if} \ \ 0 \leq t \leq \to\ - \
\eta( \epsilon); \\
T(t)(X_1^{\epsilon}  +
(t  -  \to)(X_2^{\epsilon}   - Y^{\epsilon})) & {\rm if} \
\to  - \eta( \epsilon) < t < \to  +  \eta( \epsilon); \\
X(t)   & {\rm if}\  \to -  \eta( \epsilon ) \leq t \leq b .
\end{array}
\right. \]
\end{small}
Since  $X^{\epsilon}=T(t)(Y(t)),$ where  $Y(t)$ is a
piecewise differentiable application,
except at the points $t=\to +  (\eta (\epsilon)$ and
$t=\to - \eta(\epsilon),$ we have
\[
I^{N}(X^{\epsilon},X^{\epsilon}) \geq I^{N}(T(t)(u),T(t)(u)) .
\]
On the other hand, the Focal Index of $X$ is given by
\begin{small}
$$\begin{array}{lcl}
I(X,X)&=& I(X^{\epsilon},X^{\epsilon}) \\
      &-& \int_{ \to - \eta(\epsilon) }^{ \to + \eta(\epsilon) }
          \langle \fis{X}^{\epsilon}(t),\fis{X}^{\epsilon}(t) \rangle\ - \
 \langle R({X}^{\epsilon}(t), \fis c (t))\fis c (t),{X}^{\epsilon}(t) \rangle dt \\
      &+&  \int_{ \to - \eta ( \epsilon)}^{ \to + \eta(\epsilon) }
          \langle \fis X  (t), \fis X (t) \rangle \ - \
          \langle R(X(t), \fis c (t)) \fis c (t),X(t) \rangle dt .
\end{array}$$
\end{small}
Now, using $(1)$ and $(2)$  it is easy to check that
$$
\lim_{\epsilon \so 0} I^N ( X^{\epsilon}, X^{\epsilon} )= I^N (X,X)
\geq I^N (J,J),
$$
where $J(t)=T(t)(u).$ This proves the Focal Index Lemma if there
is only one pathological point. However, one can generalize easily the above
proof to a finite number of pathological points.
\end{proof}
\begin{cor} \label{compara}
Let $(M,\langle \cdot, \cdot \rangle )$ 
be a weak Riemannian Hilbert manifold and let
$S$ and $\Sigma$ be two submanifolds
of codimension $1.$ Assume that there exists
$p \in S \cap \Sigma$ such that  $T_p \Sigma=T_p S$ is a closed subspace
of $(\tpm , \tau).$ We
denote by $N$ and $\ov{N}$ the normal vector fields
to $S$ and $\Sigma$ respectively. Suppose also that
\[
\langle \nb_{X}{N},X \rangle  < \langle \nb_{X}{\ov{N}},X \rangle,
\]
for every  $X \in T_p \Sigma =T_p S$. Then, if the Jacobi flow $T$ of
$S$ is invertible in  $(0,b)$, then the Jacobi flow of
$\Sigma$ must be injective in $(0,b).$ Moreover, if $M$ is a Riemannian
Hilbert manifold, assuming 
$A-\ov{A}$ is invertible, where
$A$ and $\ov{A}$ are the Weingarten operators at
$p$ of $S$ and $\Sigma$ respectively, 
then the Jacobi flow of $\Sigma$ is also
invertible in $(0,b).$
\end{cor}
\begin{proof} Let $s \in (0,b)$ and  let $Y(t)$ be a $\Sigma$-Jacobi field.
Since $T(t)$ is invertible in $(0,b)$, there exists a piecewise
differentiable application
$X:[0,s] \lra T_p M$  with
$X(0) \in T_p S$ such that
$Y(t)=T(t)(X(t)).$ Hence
$$
\begin{array}{rcl}
Y(0) &=& T(0)(X(0)) \\
\fis Y (0)&=& T'(0)(X(0))\ + \ T(0)(\fis X (0)) \\
(-\ov{A}(Y(0)), P_n ( \fis{Y} (0))) &=& (-A (X(0))\ +\ P (\fis X (0)),0)
\end{array}
$$
where $P_n$ is the projection of $\tpm$ onto  $T_p^{\perp} \Sigma.$
Therefore, $Y(0)=X(0)$ and the tangent component
of $\fis X (0)$ is given by $(A\ -\ \ov{A})(X(0)).$
Then,
\begin{small}
$$
\begin{array}{lcl}
\langle Y(s),\nbp{s}{Y} (s) \rangle & = & I^S (Y,Y) \\
& = & \langle (A -  \ov{A})(X(0)), X(0) \rangle +
                        \int_{0}^{s} \pr T(t)(\fis X (t)) \pr^2 dt  \\
                     & > & 0.
\end{array}
$$
\end{small}
In particular, the Jacobi flow of $\Sigma$ is injective in $(0,b).$
If $A- \ov{A}$ is invertible, then
$ \frac{d}{ds} \langle Y(s),Y(s) \rangle
\geq \pr ( A - \ov{A})^{\frac{1}{2}} \pr^{-1} \pr Y(0) \pr^2,$
so $\pr T(s)(w) \pr^2 \geq C \pr w \pr^2.$ In the Riemannian
context,  one can prove that the
image of the Jacobi flow relative to $\Sigma$ is a closed subspace for
every $s \in (0,b)$ so the Jacobi flow must be invertible in $(0,b)$
since $(Im T(s))^{\perp}= Ker T^* (s)=0,$ by Proposition \ref{p1}. 
\end{proof}
\begin{thm} \label{pippo}
Let $(M,\langle \cdot ,\cdot \rangle)$, $(N,\langle \cdot ,\cdot \rangle^*)$ be
weak Riemannian Hilbert manifolds with a Levi-Civita connections,
modeled on $\H_1$ and $\H_2$ respectively, with $\H_1$ isometric
to a closed  subspace of $\H_2$. We denote by $K^M(X,Y),$
respectively $K^N (Z,W)$ the sectional curvature of $M$ relative
the plane generated by $X,Y,$ respectively the sectional curvature
of $N$ relative the plane generated by $Z,W.$ Let
\[
c:[0,a] \lra M,\ \ c^*:[0,a] \lra N
\]
be geodesics of equal length. Suppose  that
that for every  $t\in [0,a]$
and for every  $X \in T_{c(t)}M$, $X_o \in T_{c^*(t)}N$ we have
\[
K^N (X_o,\fis c^*(t)) \geq K^M (X,\fis c(t)).
\]
Hence we have:
\begin{enumerate}
\item {\bf (Rauch)}
assume that $c^*$ has at most a
finite number of pathological points of $c^* (0)$ along $c^*$.
Let  $J$ and  $J^*$ be Jacobi fields along $c$ and $c^*$ such that
$J(0)$ and $J^*(0)$ are tangent to $c$ and $c^*$ respectively and
\begin{itemize}
\item $\pr J(0) \pr = \pr J^*(0) \pr^*$;
\item $\langle \fis c(0),\nbp {t}{J(0)} \rangle =
\langle \fis c^*(0),\nbp {t}{J^*(0)} \rangle^* $;
\item $\pr \nbp {t}{J(0)} \pr =\pr \nbp {t}{J^*(0)} \pr^*$.
\end{itemize}
Then, for every  $t\in [0,a]$
\[
\pr J(t) \pr \geq \pr J^*(t) \pr^*;
\]
\item {\bf (Berger)}
assume $c^*$ has at most a finite number of phatological
focal points of $c^* (0)$ along $c^*$, with respect the 
submanifold $N$ defined by
$\fis{c}^* (0).$ Let $J$ and
$J^*$ be Jacobi fields along $c$ and $c^*$
satisfying  $\nbp{t}{J(0)}$  and $\nbp{t}{J^* (0)}$ are
tangent to $\fis{c}(0)$ and
$\fis{c}^* (0)$ and
\begin{itemize}
\item $\pr \nbp{t}{J(0)} \pr = \pr \nbp{t}{J^* (0)} \pr^*$,
\item $\langle \fis{c}(0),J(0) \rangle = \langle
\fis{c}^* (0),J^*(0) \rangle^* ,\ \pr J(0) \pr = \pr J^*(0) \pr^*$.
\end{itemize}
Then
\[
\pr J(t) \pr \geq \pr J^* (t) \pr^*,
\]
for every $t \in [0,a]$.
\item The index of $D^2 E(c^* ) $ is greater than the index $+$ the nullity
of the index form $\D.$
\end{enumerate}
\end{thm}
\begin{proof}
We shall proof briefly the Rauch Theorem, since one may proof  likewise 
the Berger Theorem,   and we shall discuss the  third item.

%
One may note that we can assume that the Jacobi fields satisfy
\[
\pr J(0) \pr=\langle \fis c(0),\nbp {t}{J(0)} \rangle =  \pr J^*(0) \pr^*= 
\langle \fis{c}^*(0),\nbp{t}{J^*(0)}  \rangle^*=0,
\]
since the first and the second conditions imply  
$\langle J (t), \fis c (t) \rangle = \langle J^* (t), \fis{c}^*(t) \rangle^*$,
which means that the norm of the component of $J$ along $\fis c$ 
is equal to the  norm of the component of $J^*$ along $\fis{c}^*$ . 
We note also, by assumption, that $J^*(t) \neq 0$ for every $t \in (0,a]$. Let 
$\to \in (0,a]$ and let $F$ be an isometry which satisfies
$$
\begin{array}{lcl}
&F: T_{c(0)}M \lra T_{c^*(0)}N &  \\
&F(\fis c(0))=\fis {c}^*(0) &  \\
&F(\pt{\to}{0}(J(\to)))=\ptn{\to}{0}(J^*(\to))
\frac{\pr J(\to) \pr}{\pr J^*(\to) \pr^* }
\end{array}
$$
where $\ptn{s}{t}$ is the parallel transport from $c^* (s)$ to $c^* (t)$ 
along $c^*,$
and we consider the following curve of isometries
\begin{eqnarray}
& i_t: T_{c(t)}M \lra T_{c^*(t)}N \nonumber \\
& i_t=\ptn {0}{t} \circ F \circ \pt {t}{0}, \nonumber
\end{eqnarray}
$0 < t \leq t_o$. Let $W(t)=i_t(J(t))$. Put $c_o= c_{|[0, t_o]}$ 
and $c^*_o = c^*_{| [0,t_o]}.$ Then
\begin{small}
$$\begin{array}{lcl}
D^2 E(c^*_o)(W,W) &=& \int_0^{\to} {\pr \nbp{t}{W(t)} \pr^*}^2 \
                  - \ \langle R^N (\fis c^*(t),W(t))\fis c^*(t),W(t) 
\rangle^* dt  \\
              &\leq & \int_0^{\to} \pr \nbp {t}{J(t)} \pr^2
              - \langle R^M (\fis c(t),J(t))\fis c(t),J(t) \rangle dt  \\
            &=& D^2 E (c_o) (J,J) .
\end{array}$$
\end{small}
In particular we have
\begin{eqnarray}
\frac{1}{2}\d1{}{t} \mid_{t=\to} \langle J(t),J(t) \rangle &=& 
\langle J(\to),\nbp{t}{J(\to)} \rangle
 \nonumber  \\
         &=& D^2 E(c_o) (J,J) \nonumber \\
  &\geq& D^2 E(c^*_o  )(W,W) \nonumber \\
  &\geq& D^2 E(c^*_o ) (J^* \frac{\pr J(\to) \pr}{\pr J^*(\to) \pr^*}
   ,J^* \frac{\pr J(\to) \pr}{\pr J^*(\to) \pr^*}),   \nonumber
\end{eqnarray}
where the last inequality is given  by the Focal Index Lemma.
Since $J^*(t)$ is a Jacobi field we have 
$D^2 E(c^*_o )=\langle J^*(\to),\nbp{t}{J^*(\to)} \rangle^*$, so
\begin{eqnarray}
\d1{}{t} \mid_{t=\to} \langle J(t),J(t) \rangle &\geq &
   D^2 E(c^*_o )(J^* \frac{\pr J(\to) \pr}{\pr J^*(\to) \pr^*},
   J^* \frac{\pr J(\to) \pr}{\pr J^*(\to) \pr^* })\nonumber \\
        &=&
         \d1{}{t} \mid_{t=\to} \langle J^*(t),J^*(t) \rangle^*
         \frac{\pr J(\to) \pr^2 }{{\pr J^*(\to) \pr^*}^2 } . \nonumber
\end{eqnarray}
Hence given $\epsilon >0$, for every $t \geq \epsilon$ we have
\[
\d1{}{t} \log (\pr J(t)\pr^2) \geq \d1{}{t} \log ({\pr J^*(t) \pr^* }^2).
\]
Integration over $[\epsilon,t]$, yields
\[
\frac{\pr J(t) \pr^2 }{\pr J(\epsilon) \pr^2 } \geq
\frac{{\pr J^*(t) \pr^* }^2 }{{\pr J^*(\epsilon) \pr^* }^2} .
\]
Since $\pr \nbp {t}{J(0)} \pr = \pr \nbp {t}{J^*(0)} \pr^*$ 
we get our inequality.

What does it means that
the index of $D^2 E(c^* ) $ is greater than the index $+$ the nullity
of $\D?$

Let $ i_o:T_{c(0)} M \lra T_{c^* (0)}N$ be an
isometry such that $i_o (\fis{c}(0))= \fis{c}^* (0)$.
We define, for each $t \in [0,b],$ the isometry
$$
i_t = \pt{0}{t} \circ i_o \circ \ptn{t}{0}:
T_{c(t)}M \lra T_{c^* (t)} N.
$$

Let $X$ be a vector field along $c.$ We may consider the vector field
$i(X)(t)= i_t (X(t))$ along $c^*$
and one can prove that
$$
\D (X,X) \geq D^2 E (c^*) (i(X),i(X)),
$$ 
by the assumption on the sectional curvatures.
Then if $U$ is a subspace on which $\D \leq 0$ then
$D^2 E (c^*)_{i (U)} \leq 0.$
\end{proof}

\section{The free loop space of a finite dimensional Riemannian manifold}
Let $(M^n, \langle \cdot, \cdot \rangle )$ be a compact Riemannian manifold of
dimension $n.$ We recall that 
$\Omega (M^n)=H^1 (S^1, M^n)$ is the set of maps of Sobolev
class $H^1$ from $S^1$ into $M^n.$ It can be given the structure of
an infinite dimensional Hilbert manifold and the tangent space
$T_{\sigma} \Omega (M^n)$ at a point $\sigma \in \Omega (M^n)$
consist of periodic $H^1$ vector fields along $\gamma.$ One defines the
$L^2$  weak Riemannian structure on $\Omega (M^n)$ by setting
$$
\langle X,Y \rangle (\sigma) = \int_{S^1} \langle X(t), Y(t) \rangle dt
$$
where $X,Y \in T_{\sigma} \Omega(M^n).$
It is well known, see  \cite{df} or \cite{Mi3}, that the $L^2$ metric has a
Levi-Civita connection which is determined pointwise by the Levi-Civita
connection of $M^n.$ Moreover, the $L^2$ curvature $R$ is
given pointwise by the
tensor curvature of $M^n,$ so the sectional curvature is given by
$$
K(X,Y)= \int_{S^1}  K^{M^n} (X(t), Y(t))  dt.
$$
If $M^n$
has positive sectional curvature, i.e. $K^{M^n} \geq K_o > 0,$ then
$\Omega (M)$ has positive sectional curvature
since $K^{\Omega(M^n)} \geq K_o 2 \pi=K_1$. In particular there exist
at least a conjugate point along any geodesic of
length greater than  $\frac{\pi}{\sqrt{K_1}}$ and its index is infinite.
Indeed, let $\gamma:[0,1] \lra \Omega(M^n )$ be a geodesic with length
$l> \frac{\pi}{\sqrt{K_1}}.$ Let $v \in T_{\gamma (0)} \Omega(M^n )$
be a unit vector such that $ \langle v,w \rangle (\gamma (0) ) =0.$
Let $W(t)=\sin (t \pi) V(t),$ where $V$ is the parallel transport
along $\gamma$ of $v.$ One can verify that
$D^2 E (\gamma) (W,W) <0,$ so, by
Focal Index Lemma, we have at least a
singularity of the exponential map and it cannot be an isolated pathological
point. The fact that the index is infinite follows comparing $\Omega (M^n)$ with 
the manifold
$$
 S_{\frac{1}{\sqrt{K_1}}} := \{ (x,Y) \in \R \times T_{\gamma(0)} 
\Omega (M^n) : \ x^2 +  \langle Y, Y \rangle (\gamma (0))= 
\frac{1}{\sqrt{K_1}} \}
$$
which is a weak Riemannian Hilbert manifold of a constant sectional
curvature $K_1.$ 
Note that the same argument works if we consider $N$ the sub\-ma\-ni\-fold 
defined
by $\fis{\gamma} (0).$ Indeed, one 
verifies that $D^2 E (\gamma) (W,W) <0,$
where $W(t)= \cos (t \frac{\pi}{2}) V(t).$
This fact is
in contrast with the Riemannian point of view, i.e. $\Omega (M^n)$
endowed by the $H^1$ metric, since Misiolek  proved \cite{Mi} that the
exponential map is a non-linear Fredholm map 
and any geodesic of finite length has finite index.

Suppose now $M^n =G$ is a non-abelian compact Lie group. In this case
we get a simple expression for the Levi-Civita connection
$$
\nabla_X Y = \frac{1}{2}[X,Y]
$$
and therefore, for the curvature tensor
$$
R(X,Y)Z= -\frac{1}{4}[[X,Y],Z].
$$
Consequently, any one-parameter subgroup of $\Omega (G)$ is a geodesic of
the $L^2$ metric and the exponential map is defined on the whole 
tangent space.
Moreover, if $X,Y$ and $Z$ are parallel vector fields along a geodesic $c$
then $R(X,Y)Z$ is parallel along $c$ as well
(see \cite{mi} and \cite{mi2}).

It is well known, see \cite{He}, that
Lie $(G)$= $\mathfrak{z} \oplus \mathfrak{g}_s,$  where $\mathfrak{g}_s$
is the maximal semisimple ideal of Lie$(G)$ and $\mathfrak{z}$ is the Lie
algebra of the center of $G.$ Since
$\mathfrak{g}_s$ is semisimple, it has a subalgebra
$\mathfrak{h}_{\alpha}$ isomorphic to $\mathfrak{su} (2).$
We denote by $A_{\alpha},$
$B_{\alpha}$ and  $C_{\alpha}$ the standard generators of
$\mathfrak{su} (2).$ Then
$$
[A_{\alpha}, B_{\alpha}]=2C_{\alpha},
[C_{\alpha}, A_{\alpha}]=2B_{\alpha},
[C_{\alpha}, B_{\alpha}]=-2A_{\alpha}.
$$
Let $c$ be the one-parameter subgroup of $\Omega (G)$ generated by
$\frac{1}{\sqrt{2 \pi}}{B}_{\alpha}.$ Now, as in \cite{Mi2}
p.2480--2481, one can prove that the vector fields $Y_k (t) = \sin
(\frac{t}{\sqrt{2 \pi}} ) \pt{0}{t} (w_k ),$ where $w_k (x) =
\frac{1}{\sqrt{\pi}} \sin k x A_{\alpha}$ is an eigenvector of $R
(\cdot , \frac{1}{\sqrt{2 \pi}}B_{\alpha}) \frac{1}{2 \pi
}B_{\alpha}$ with the eigenvalue $\lambda= \frac{1}{\sqrt{\pi}},$
are linearly independent Jacobi fields along $c.$ Then the kernel
of $\mathrm{ d( \e e)_{\pi \sqrt{2 \pi} B_{\alpha}}}$ is infinite
dimensional. In particularly $\Omega (G)$ has at least a
monoconjugate point along $c$ and the exponential map cannot be
Fredholm. Moreover, the vectors fields $L_k (t) =  \cos
(\frac{t}{\sqrt{2 \pi}} ) \pt{0}{t} (w_k ),$ are $N-$Jacobi fields
along $c,$ where $N$ is the submanifold defined by $\fis{c} (0),$
so  the kernel $\mathrm{ d( \Ex )_{ \frac{\pi \sqrt{2 \pi}}{2}
B_{\alpha}}}$ is infinite dimensional. Hence,  there exists a
monofocal point along $c$ and the application $\Ex: T^{\perp} N
\lra M,$ which is well-defined in this case, fails to be Fredholm.

The author thank Professors Francesco Mercuri and  Paolo Piccione for many 
useful conversations on the subject, as well the referee for helpful 
suggestions. I thank Professors Gerard Misiolek and Christopher J.
Atkin  for useful remarks. I also thank Professor Alessandro Ghigi 
for his support and for his continuing interest in my work.

\end{document}